\documentclass[english]{amsart}

\usepackage[all,cmtip]{xy}
\usepackage{hyperref}
\usepackage{blindtext}

\usepackage{geometry}
\usepackage{amsmath,amssymb,amscd,amsfonts}

\usepackage{babel}
\usepackage{amstext}
\usepackage{amsmath}
\usepackage{amsfonts}
\usepackage{latexsym}
\usepackage{ifthen}

\usepackage[all,cmtip]{xy}
\xyoption{all}
\usepackage{xcolor}
\newcommand\chr[1]{{\textcolor{blue}{#1}}}

\usepackage{enumitem}
\setlist[enumerate]{label=(\thethm.\arabic*), before={\setcounter{enumi}{\value{equation}}}, after={\setcounter{equation}{\value{enumi}}}}

\newcommand{\R}{\mathbb{R}}
\newcommand{\C}{\mathbb{C}}
\newcommand{\Q}{\mathbb{Q}}

\newcommand{\ddbar}{\partial\bar{\partial}}

\newcommand{\cC}{\mathcal{C}}

\newcommand{\cI}{\mathcal{I}}

\renewcommand{\O}{\mathcal{O}}

\newcommand{\ep}{\varepsilon}
\renewcommand{\epsilon}{\varepsilon}

\newcommand{\ol}{\overline}

\renewcommand{\leq}{\leqslant}
\renewcommand{\geq}{\geqslant}
\newcommand{\Ricci}{\mathrm{Ricci}}
\newcommand{\codim}{\mathrm{codim}}

\newcommand{\Supp}{\mathrm {Supp}}
\newcommand{\tr}{\mathrm{tr}}

\newcommand{\Psh}{\mathrm{Psh}}

\newcommand{\dbar}{\bar \partial}

\newtheorem{thm}{Théorème}[section]
\newtheorem{lemma}[thm]{Lemme}

\newtheorem{conjecture}[thm]{Conjecture}
\newtheorem{defn}[thm]{Definition}
\newtheorem{cor}[thm]{Corollaire}

\newtheorem{remark}[thm]{Remarque}

\numberwithin{equation}{thm}

\title{Fonctions plurisousharmoniques et g\'eom\'etrie complexe: \\ sur quelques r\'esultats de J.-P. Demailly }
\author{Mihai PAUN}

\address{Universit\"at Bayreuth, Mathematisches Institut, Lehrstuhl Mathematik VIII, Universit\"atsstrasse 30, D-95447, Bayreuth, Germany}
\email{mihai.paun@uni-bayreuth.de}

\begin{document} 

\maketitle

%\begin{dedication}
%\vspace{-1.5cm}
%{To Bo Berndtsson,\\ On the Occasion of His $70^{\rm th}$ Birthday}
%\end{dedication} 
\vskip 1cm

\section{Introduction}

La th\'eorie des fonctions plurisousharmoniques (``psh'' en abr\'eg\'e) joue un rôle déterminant dans l'\oe uvre math\'ematique de J.-P. Demailly. Le but de cet article est de pr\'esenter certains de ses résultats où ceci est particulièrement
perceptible. Les travaux qu'on va survoler par la suite sont, en quelque sorte, représentatifs pour un des aspects importants de la stratégie de recherche de Jean-Pierre, i.e. implémenter systématiquement des méthodes analytiques en géométrie algébrique. Nous allons essayer d'illustrer la richesse extraordinaire de ses idées en mentionnant à la fin de chaque section quelques références qui montrent l'impact de ses résultats sur le développement du domaine.
\medskip

Remarquons avant tout que Jean-Pierre a été très tôt initié aux méthodes analytiques, honneur à son directeur de thèse, M. Henri Skoda ! Comme preuve, citons un des ses premier résultats, publié dans \cite{Dem79}. \footnote{alors qu'il avait que vingt-deux ans...}.

\begin{thm}\label{dem1} Soit $S\subset \C^2$ la surface d'\'equation
 \[e^x+ e^y= 1.\]
Alors $S$ est de Liouville, i.e. toute fonction holomorphe born\'ee sur $S$ est constante. 
\end{thm}
En fait, le théorème principal dans \cite{Dem79} est plus général : \emph{si $f$ est une fonction holomorphe d\'efinie sur $S$, telle qu'il existe un entier $n$ et une constante $C> 0$ avec $$|f(z)|\leq C(1+ |z|)^n,\qquad z=(x,y)\in S$$ alors $f$ est la restriction \`a $S$ d'un polynôme $P(x,y)$ de degré total au plus $n$.} Le preuve repose essentiellement sur les estimées $L^2$ de H\"ormander, combinées avec le fait que $S$ est ``proche`` de la surface définie par l'équation $e^x+ e^y= 0$ qui consiste en une réunion de droites. \footnote{\`A titre d'anecdote, je me suis retrouvé aux alentours de 2005/06 à Oberwolfach en compagnie de F. Bogomolov et M. McQuillan, qui m'ont demandé si j'étais familier avec le résultat \ref{dem1}... car ils envisageaient une preuve du théorème de Fermat en l'utilisant !}
\medskip

Ce qui suit est divisé en trois parties. Nous avons choisi ``d'attaquer'' par un article de Jean-Pierre (accidentellement, publié dans les proceedings d'une conférence à Bayreuth, cf. \cite{DemBT}) où on voit déjà germer les notions et résultats qui vont jouer un rôle crucial dans la suite. Ensuite, nous allons expliquer quelques idées
importantes dans l'article \cite{Dem93}, qui représente une des contributions fondamentales de Jean-Pierre. Finalement, nous allons discuter les étapes principales de la preuve du théorème principal dans \cite{DP04}. i.e. la caractérisation numérique du cône de K\"ahler d'une variété kählérienne compacte.
\bigskip

\section{M\'etriques singuli\`eres et g\'eom\'etrie : le dictionnaire de J.-P. Demailly}

Soit $X$ une variété complexe compacte, et soit $L$ un fibré holomorphe en droites. Un invariant très important de $L$ est sa classe de Chern 
notée $c_1(L)$; il y a --au moins-- deux manières de l'introduire. 

Considérons d'abord le point de vue hermitien, soit donc $h_L$ une
métrique sur $L$. Cela signifie que pour toute trivialisation
\[\theta: L|_\Omega\to \Omega\times \C\]
de la restriction de $L$ sur un ouvert $\Omega\subset X$ on dispose d'une fonction poids $\varphi\in \cC^\infty(\Omega, \R)$ telle que l'égalité
\begin{equation}\label{eq1}
|\xi|_{h_L}^2= |\theta(\xi)|^2e^{-\varphi(x)}\end{equation}
soit vérifiée pour tout $x\in \Omega$ et pour tout $\xi \in L_x$. Autrement dit,
la fonction $\varphi$ mesure l'écart entre la métrique plate sur le fibré trivial $\Omega\times \C$ et $h_L|_\Omega$.

Etant donné un fibré hermitien $(L, h_L)$, la forme de courbure correspondante sera notée par $i\Theta(L, h_L)\in \cC^{\infty}_{1,1}(X, \R)$. Ladite forme s'exprime localement par la formule
\[i\Theta(L, h_L)|_\Omega= i\ddbar \varphi,\]
et la classe de cohomologie induite par la forme globale $\displaystyle \frac{i}{2\pi}\Theta(L, h_L)$ représente la première classe de Chern de $L$.
\smallskip

Par ailleurs, soit $s$ une section meromorphe de $L$. Il se trouve que la classe de cohomologie définie par le diviseur $s=0$ vaut également $c_1(L)$.
\smallskip

\noindent Dans ce contexte, la généralisation suivante de la notion de métrique est naturelle, car elle unifie les deux points de vue mentionnés précédemment.

\begin{defn}\label{def1}
  Une métrique singuli\`ere sur $L$ est définie par la formule \eqref{eq1}
  via des fonctions poids
  $\varphi\in L^1_{\rm loc}(\Omega)$.
\end{defn}

Autrement dit, la seule différence par rapport au cadre usuel est la régularité des fonctions poids qui définissent $h_L$. Remarquons qu'on
peux aussi définir la notion de ``courbure'' d'un fibré muni d'une métrique singulière. En général, ce n'est plus une forme différentielle, mais plut\^ot un courant fermé \footnote{Je ne vais pas définir de façon formelle cette notion, voir e.g. \cite{bookJP}, chapitre trois.}
$T$ défini comme suit.
Soit $\beta\in \cC^{\infty}_{n-1, n-1}(X)$ une forme différentielle de type $(n-1, n-1)$ sur $X$, et supposons que son support est contenu dans un ouvert de trivialisation $\Omega\subset X$. Alors $T$ agit sur $\beta$ selon la formule 
\begin{equation}\label{eq2}
\langle T, \beta \rangle:= \frac{1}{2\pi}\int_X \varphi i \ddbar \beta,
\end{equation}
ce qui a un sens, compte tenu des propriétés de $\varphi$ dans la définition \ref{def1}. On étend l'égalité \eqref{eq2} pour toute forme 
$\beta$ (à support compact lorsque $X$ n'est pas compacte) par une partition de l'unité. Le théorème de Stokes montre que $T$ est bien défini, i.e. il ne dépend pas de la manière dont on choisit les trivialisations de $L$.
\smallskip

\subsection{Quelques exemples}
\begin{enumerate}

\item[(1)] Toute métrique hermitienne sur $L$ est en particulier une métrique singulière, dont le courant de courbure coïncide avec $i\Theta(L, h_L)$,
i.e. 
\[\langle T, \beta \rangle= \int_X i\Theta(L, h_L)\wedge \beta. \]

\item[(2)] Soit $D= \sum e_i D_i$ un diviseur sur $X$, ou les $D_j$ sont des hypersurfaces irréductibles et les coefficients $e_i\in \mathbb Z$ sont des entiers. On utilise la même notation $L_D$ pour designer le faisceau $\mathcal O(D)$ dont les sections locales sont les fonctions méromorphes $f$ telles que $(f)+ D\geq 0$ et le fibré en droites 
associé. Il est aisé de construire les trivialisations locales de $L_D$:
\[\theta(f):= f\prod g_j^{e_j}\]  
où localement on a $D_j= (g_j=0)$. 

On introduit la métrique singulière $\displaystyle \Vert f\Vert:= |f|$ ; par rapport à la trivialisation ci-dessus, son poids vaut
$\displaystyle \varphi= \sum e_j\log|g_j|^2$
et par la formule de Poincaré-Lelong, le courant de courbure s'exprime comme suit
\[\langle T, \beta \rangle= \sum e_j\int_{D_{j, \rm reg}}\!\!\!\!\!\beta,\]
i.e. $\displaystyle T= \sum e_j[D_j]$, le courant d'intégration sur $D$.

\item[(3)] Soient $s_1,\dots s_N$ sections holomorphes d'une puissance tensorielle $L^k$ de $L$, où $k\geq 1$ est un entier positif. On peux leur associer une métrique singulière $h_L$ sur $L$ par la formule
\[|\xi|_{h_L}^2:= \frac{|\theta(\xi)|^2}{(\sum |f_j|^2)^{\frac{1}{k}}},\]
où les fonctions $f_j$ correspondent aux sections $s_j$ via les trivialisations $\theta^{\otimes k}$ de $L^k$. Le poids (local) de $h_L$ est donc 
$\displaystyle \varphi= \frac{1}{k}\log(\sum |f_j|^2).$ 
\end{enumerate}

\medskip

\noindent Considérons un domaine $\Omega\subset \mathbb C^n$ dans l'espace euclidien. Parmi les fonctions dans $L^1_{\rm loc}(\Omega)$ on distingue la classe suivante. 

\begin{defn}\label{def2} Une fonction $\varphi: \Omega\to [-\infty, \infty[$ est dite plurisousharmonique ("psh" en abrégé) si les propriétés suivantes sont satisfaites :
\begin{itemize}

\item $\varphi$ est semicontinue supérieurement, et non-identiquement $-\infty$ sur chaque composante de $\Omega$ ;

\item la restriction de $\varphi$ à chaque droite complexe vérifie l'inégalité de la moyenne, i.e. $\forall z_0\in \Omega$ et quelque soit
$\xi\in \mathbb C^n$ tel que $z_0+ e^{i\theta}\xi\in \Omega$ pour tout $\theta\in [0, 2\pi]$ on a 
\[\varphi(z_0)\leq \frac{1}{2\pi}\int_0^{2\pi}\varphi(z_0+ e^{i\theta}\xi)d\theta.\]
\end{itemize}
L'ensemble des fonctions psh définies sur $\Omega$ sera noté $\Psh(\Omega)$.
\end{defn}
\medskip

\noindent Voici quelques propriétés de ces fonctions,\footnote{Lors d'un cours d'analyse complexe par Jean-Pierre, il nous a présenté certaines de ces propriétés en faisant le commentaire \emph{bon, vous voyez, ce sont presque des fonctions de classe} $\mathcal C^\infty$"...} qui font qu'elles seront particulièrement intéressantes.

\begin{thm} Soit $\Omega\subset \mathbb C^n$ un domaine de $\mathbb C^n$, et soient $\varphi, \varphi_1, \varphi_2\in \Psh(\Omega)$
ainsi que $\psi\in L^1_{\rm loc}(\Omega)$, semi-continue supérieurement. Les affirmations suivantes sont vraies.
\begin{enumerate}
\smallskip

\item[\rm (1)] Soit $(f_j, r_j)_{j=1,\dots, m}\subset \mathcal O(\Omega)\times \R_+$ un sous ensemble fini. Les fonctions 
\[\varphi:= \log\big(\sum |f_j|^{2r_j}\big), \qquad \psi:= -\log\log\frac{1}{\sum |f_j|^{2r_j}}\]
sont psh.
\smallskip

\item[\rm (2)] On a $\varphi\in L^p_{\rm loc}(\Omega)$ pour tout $p> 0$, et $\max(\varphi_1, \varphi_2)\in \Psh(\Omega)$.
\smallskip

\item[\rm (3)] Le courant $\displaystyle T:= i\ddbar \varphi$ est positif.
\smallskip

\item[\rm (4)]  Si $i\ddbar \psi\geq 0$ au sens des courants sur $\Omega$, alors $\psi\in \Psh(\Omega)$.
\smallskip

\item[\rm (5)] Pour tout $x\in \Omega$ il existe $\gamma> 0$ tel que $\displaystyle \int_{(\Omega_, x)}e^{-\gamma \varphi}d\lambda< \infty$.
\smallskip

\item[\rm (6)] Pour tout $\ep> 0$, le gradient de $\varphi$ au sens des distributions est dans $L^{2-\ep}_{\rm loc}$.

\end{enumerate}
\end{thm}

\noindent Les propriétés (1)-(4) sont élémentaires, alors que (5) et (6) sont un peu plus subtiles (cf. e.g. \cite{bookJP} et les références dedans).

\medskip

\noindent Pour les lecteurs qui ne sont pas (encore) familiers avec la théorie des courants, la propriété (3) ci-dessus 
peux se comprendre comme suit : il existe un ensemble de mesures $(\mu_{j\ol k})_{1\leq j,k\leq n}$ telles que
 \begin{equation}\label{eq5}\overline{\mu_{j\ol k}}= \mu_{k\ol j}, \qquad \sum_{j, k}\int_\Omega \theta_j\overline \theta_k\mu_{j\ol k}\geq 0\end{equation}
ainsi que 
\begin{equation}\label{eq4}
\langle T, \beta\rangle =  \int_\Omega \sum_{j, k} i\mu_{j\ol k}dz_j\wedge dz_{\ol k}\wedge\beta \end{equation}
pour toute forme $\beta$ de bidegré $(n-1, n-1)$ et tout ensemble de fonctions régulières $(\theta_k)$ à support compact sur $\Omega$.
On rappelle que par définition on a $\displaystyle \langle T, \beta\rangle = \int_\Omega\varphi i\ddbar\beta$.
\smallskip

Dans le cas des diviseurs, on dispose de la notion de "multiplicité" ; l'analogue pour les fonctions psh sont les nombres de Lelong,
définis comme suit.

\begin{defn}\label{lolo}
Soient $\varphi\in \Psh(\Omega)$, et $x\in \Omega$ un point arbitraire. On appelle le nombre de Lelong de $\varphi$ en $x$ la limite inf suivante 
\[\nu(\varphi, x):= \lim\inf_{z\to x}\frac{\varphi(z)}{\log|z-x|}.\] 
\end{defn}

On peux calculer $\nu(\varphi, x)$ en utilisant le courant $\displaystyle T= \frac{i}{\pi}\ddbar \varphi$ associé à $\varphi$ ; la formule est la suivante
\begin{equation}\label{eq3}
\nu(\varphi, x)= \lim_{r\to 0_+}\frac{1}{(2\pi r^2)^{n-1}}\int_{|z-x|< r}T\wedge \omega_{\rm euc}^{n-1}=: \nu(T, x)\end{equation}
où $\omega_{\rm euc}:= i\ddbar |z|^2$. Lorsque $\varphi= \log|f|^2$, le nombre de Lelong de $\varphi$ en $x$ coïncide avec la multiplicité
du diviseur $f=0$ en $x$, cf. \cite{bookJP} et compte tenu de la formule \eqref{eq3}, cela offre une très jolie interprétation du dernier. 
\smallskip

Retournons à présent au cadre global d'une variété complexe compacte $(X, \omega)$ munie d'une métrique hermitienne. Via les coordonnées locales de $X$ la notion de fonction psh sur $X$ a encore un sens. Mais on voit aisément que l'espace $\Psh(X)$ est réduit aux fonctions constantes (par l'inégalité de la moyenne). Du coup, on peux considérer que 
l'analogue global des fonctions psh sont les fonctions $\phi: X\to [-\infty, \infty[$, semicontinues supérieurement et non-identiquement $-\infty$, telle qu'il existe une constante $C>0$ de sorte qu'on ait
\[T= i\ddbar \phi+ C\omega\geq 0\]
au sens des courants sur $X$. Autrement dit, étant donné $x\in X$ et un système de coordonnées locales $(U, z)$ centré en $x$, il existe une
constante positive $C_x> 0$ telle que $\displaystyle \phi|_{U}+ C_x|z|^2\in \Psh(U).$ Une telle fonction $\phi$ sera appelée \emph{quasi-psh}. Si on préfère, on peux considérer directement 
l'ensemble de courants fermés $T$ de type $(1,1)$, tels que $T\geq -C\omega$ sur $X$. Dans e.g. \cite{Dem92} on montre l'existence d'une 
forme différentielle type $(1,1)$ disons $\gamma$, réelle et fermée, ainsi que d'une fonction quasi-psh $\phi$ telle que l'égalité suivante
\[T= \gamma+ i\ddbar\phi\] 
soit vérifiée (au sens faible). 
\smallskip

Il se trouve que l'analogie entre diviseurs et courants associés aux fonctions psh est beaucoup plus profonde, compte tenu du résultat suivant
dû à Y--T. Siu.

\begin{thm}\label{siu}\cite{Siu74}
Soit $X$ une variété complexe, et soit $T$ un courant positif fermé sur $X$ de bidegré $(1,1)$. 
\begin{enumerate}
\smallskip

\item[\rm (1)] Pour chaque $c> 0$ on définit l'ensemble
\[E_c(T):=\{x\in X : \nu(T, x)\geq c\}.\]
Alors $E_c(T)$ est un sous-ensemble analytique (fermé) de $X$. 

\item [\rm (2)] Soient $(Y_j)$ les composantes en 
codimension un contenues dans les ensembles de niveau $\bigcup_{c> 0} E_c(T)$.
Il existe un courant positif fermé $R$ tel que les composantes irréductibles de $E_c(R)$ ont codimension au moins deux dans $X$, et tel qu'on a l'égalité 
\begin{equation}\label{eq5}T= \sum \nu_j[Y_j]+ R,\end{equation}
où les $\nu_j$ sont des réels positifs. 
\end{enumerate}
\end{thm}

\noindent Le résultat de Siu est encore plus complet  -et impressionnant- que cela : il reste valable pour les courants $T$ de bidegré $(p,p)$
et dans ce cas les $Y_j$ dans \eqref{eq5} sont les ensembles de codimension $p$ de $\displaystyle \bigcup_{c> 0} E_c(T)$, et les composantes de 
$E_c(R)$ ont codimension au moins $p+1$.

\subsection{Quelques résultats} L'énoncé suivant, cf. \cite{Dem82}, représente une version "globale" très poussée des estimées $L^2$ de Hörmander, cf. \cite{Hor}. C'est un outil  
fondamental en géométrie complexe.
 Soit $L\to X$ un fibré en droites, muni d'une métrique singulière $h_L$ dont les les fonctions poids sont quasi-psh.  
Supposons que les poids ${\varphi_L}$ de $h_L$ sont logarithmiques, 
i.e. localement on a ${\varphi_L}= c\log(\sum |f_\alpha|^2)+ \mathcal O(1)$, ou $(f_\alpha)$ sont holomorphes et $c\geq 0$. Alors la décomposition 
\eqref{eq5} du courant de courbure est particulièrement simple
\[i\Theta(L, h_L)= \sum_{j=1}^N \nu_j[Y_j]+ \Gamma,\]
où les hypersurfaces $Y_j$ coïncident localement avec les composantes de codimension un en $X$ de l'ensemble analytique $f_\alpha = 0$, les coefficients $\nu_j\geq 0$ et de plus
\begin{equation}\label{eq7}\Gamma|_{\Omega}= i\sum \theta_{j\ol k}dz_j\wedge d\ol z_k\end{equation}
 où les fonctions $\theta_{j\ol k}$ sont dans $L^1_{\rm loc}$ et $\Omega\subset X$ est un ouvert de coordonnées.
 \smallskip
 
 \noindent Dans ce contexte on a le résultat suivant. 
 
 \begin{thm}\label{DemL2} \cite{Dem82} Soit $X$ une variété kählérienne compacte, et soit $(L, h_L)$ un fibré en droites muni d'une métrique 
 singulière dont les poids ont des fonctions psh à singularités logarithmiques. Soit $v$ une $(n, 1)$--forme $L^2$ à 
 valeurs dans $L$. On suppose que les hypothèses suivantes 
\begin{equation}\label{eq6}\dbar v=0, \qquad \int_X|v|_{\Gamma, h_L}^2dV<\infty\end{equation}
sont satisfaites. Alors il existe une $(n, 0)$--forme $u$ à valeurs dans $L$ telle que
\begin{equation}\label{eq8}\dbar u=v, \qquad \int_X|u|_{h_L}^2dV\leq \int_X|v|_{\Theta\chr{\Gamma ???}, h_L}^2dV\end{equation}
\end{thm}
La notation dans \eqref{eq6} signifie qu'on mesure la forme $v$ en utilisant la partie absolument continue $\Gamma$ du courant de courbure,
i.e. $|v|_{\Gamma, h_L}^2dV$ est une forme de degré maximal qui s'écrit localement
\[|v|_{\Gamma, h_L}^2dV|_\Omega= \sum_{k,m} v_{k}\ol{v_{m}}\theta^{\ol k m}e^{-\varphi}d\lambda\] 
si $\displaystyle v|_\Omega= \sum v_k dz\wedge d\ol z_k$ ; on note $(\theta^{\ol k m})_{m,k}$ l'inverse de la matrice des coefficients de $\Gamma$ 
dans \eqref{eq7}. Finalement, même si l'hypothèse "$X$ kählérienne" est nécessaire dans la preuve, la métrique de $X$ n'intervient pas dans les estimées 
\eqref{eq8}.

\begin{remark}{\rm 
Le théorème \ref{DemL2} s'inscrit naturellement dans la pléthore des résultats concernant l'annulation de certains groupes de cohomologie 
sous hypothèses de positivité. Ce qui fait la différence ici c'est que le courant de courbure est seulement supposé semi-positif,
e.g. même si $h_L$ est non-singulière, le noyau de la courbure peut être non-trivial, pourvu que l'hypothèse d'intégrabilité dans
\eqref{eq6} est satisfaite. Ceci est crucial dans la preuve de la \emph{version ultime} du théorème d'extension d'Ohsawa-Takegoshi, cf. \cite{CDM}. Aussi, l'estimée de la norme de $u$ dans \eqref{eq8} est très précise et extrêmement utile : par exemple, le fait d'avoir la constante égale à 1 (et non pas 15) est fondamentale dans la preuve de la positivité des images directes des fibrés adjoints, cf. \cite{BoB}.  
Remarquons aussi que le théorème 4.1 dans \cite{Dem82} est plus général que l'énoncé ci-dessus : par exemple, il y a pas de restrictions sur la nature des singularités de $h_L$.}
\end{remark}
\medskip

\noindent Dans le reste de cette section nous allons présenter quelques applications du théorème \ref{DemL2}. Supposons que $(L, h_L)$ 
est un fibré en droites muni d'une métrique singulière sur une variété kählérienne compacte $(X, \omega)$ de dimension $n$, tel que
\begin{equation}\label{eq9}
T:= \frac{i}{\pi}\Theta(L, h_L)\geq \ep\omega, \qquad \nu (T, x_0)\geq n+s
\end{equation}
où $\ep > 0$ est un réel positif, $s\in \mathbb N_+$ et $x_0\in X$ est un point arbitraire. Le résultat suivant a été obtenu dans \cite{DemBT}.

\begin{thm}\label{DemGen} \cite{DemBT}
Supposons que les hypothèses \eqref{eq9} sont satisfaites, et qu'il existe un ouvert $x_0\in V\subset X$ contenant le point 
$x_0$, tel que $e^{-\varphi_L}\in L^1_{\rm loc}(V\setminus x_0)$, où $\varphi_0$ est le poids de
$h_L$ près de $x_0$. Alors l'espace des sections globales $H^0(X, K_X+ L)$ engendre les $s$-jets en $x_0$. 
\end{thm}

Afin de mieux comprendre la signification de ce résultat, esquissons la preuve. Pour commencer, l'hypothèse \eqref{eq9} concernant le nombre de Lelong en $x_0$ implique l'existence d'une constante $C> 0$ telle que l'inégalité \[e^{-\varphi_L (z)}\geq \frac{C}{|z- x_0|^{2n+2s}}\]
soit satisfaite pour tout $z$ au voisinage de $x_0$. On suppose par ailleurs que $e^{-\varphi_L}$ est localement intégrable en tout
point $x\neq x_0$ ; cela fait que si on considère $P\in \mathbb C[z]$ un polynôme arbitraire, et $\theta$ une fonction tronquante, qui vaut 1 
près de $x_0$ et dont le support est contenu dans $V$, alors la forme   
\[v:= \dbar(\theta P)\otimes e\]
satisfait les conditions \eqref{eq6}, où on note $e$ une section holomorphe non nulle en $x_0$ de $K_X+L|_V$ (quitte à restreindre $V$, on peux supposer que $K_X$ et $L$ sont triviaux sur $V$ ; ainsi, on peux trouver sans peine des sections $e$ qui n'ont pas de zeros). D'après le théorème \ref{DemL2},
on peut résoudre l'équation $\dbar u= v$, et de plus $u$ est intégrable par rapport à $h_L$. \'Etant donné que $v$ vaut zéro près de 
$x_0$, notre solution $u$ sera holomorphe au voisinage de ce point. La condition $L^2$ force $u$ s'annuler à l'ordre au moins $s+1$, du coup
le $s$-jet de la section holomorphe $\theta P\otimes e- u$ de $K_X+L$ coïncide avec $P$. Voila!
\medskip

Un fibré en droites $L\to X$ est dit \emph{ample} si l'application définie par l'espace des sections globales d'une de ses puissances tensorielles est un plongement de $X$ dans un espace projectif. Un
des théorèmes fondamentaux de K. Kodaira, cf. \cite{K3}, montre que cette propriété admet une caractérisation métrique : \emph{$L$ est ample si et seulement s'il admet une métrique hermitienne $h_L$ dont la forme de courbure est supérieure à un multiple strictement positif d'une métrique.}  

Les notions suivantes sont d'une grande importance en géométrie algébrique (et tout particulièrement dans le programme des modèles minimaux).
\begin{defn}\label{mo} Un fibré en droites $L\to X$ sur une variété projective est :
\begin{enumerate}
\smallskip

\item[\rm (a)] nef, si $\displaystyle L\cdot C= \int_{C_{\rm reg}}c_1(L)\geq 0$ pour toute courbe $C\subset X$,
\smallskip

\item[\rm (b)] pseudo-effectif, si $c_1(L)$ est limite des classes induites par des $\Q$-diviseurs effectifs, 
\smallskip

\item[\rm (c)] big,\footnote{ou "gros" en bon français, bien que cette terminologie semble avoir du mal à s'imposer...}
lorsque l'ordre de croissance de la dimension de l'espace de sections de ses puissances tensorielles est maximale, i.e. 
$h^0(X, mL)\simeq m^n$ pour $m\to\infty$. 
\end{enumerate}
\end{defn} 
\medskip

En poursuivant le théorème de Kodaira sus-mentionné, Jean-Pierre a établi dans \cite{DemBT} le "dictionnaire" suivant.

\begin{thm}\label{DemDict} Soit $(X, \omega)$ une variété projective munie d'une métrique $\omega$, et soit $L\to X$ un fibré en droites. On a les équivalences suivantes.
\begin{enumerate}
\smallskip

\item[\rm (a)] $L$ est nef si et seulement s'il admet une famille de métriques hermitiennes $(h_\ep)_{\ep> 0}$ sur $L$ telle que 
$\displaystyle i\Theta(L, h_\ep)\geq -\ep\omega$ sur $X$, pour tout $\ep> 0$.
\smallskip

\item[\rm (b)] $L$ est pseudo-effectif si et seulement s'il admet une métrique singulière $h_L$ dont le courant de courbure est 
semi-positif, i.e. $\displaystyle i\Theta(L, h_L)\geq 0$.
\smallskip

\item[\rm (c)] $L$ est big si et seulement s'il admet une métrique singulière $h_L$ dont le courant de courbure est supérieur a un multiple positif 
de $\omega$, i.e. $\displaystyle i\Theta(L, h_L)\geq \ep_0\omega$ (dans ce cas on appelle $\displaystyle i\Theta(L, h_L)$ un courant kählérien).
\end{enumerate} 
\end{thm}

\begin{remark}{\rm 
Les puissances tensorielles d'un fibré en droites pseudo-effectif n'admettent pas -en général- des sections holomorphes globales (non-identiquement nulles...). Ainsi, le point (b) dans le théorème \ref{DemDict} offre une excellente alternative, i.e. l'existence d'une métrique 
dont le courant de courbure est semi-positif : même si on ne dispose pas forcement de sections holomorphes, on a un objet sous la main pour travailler!   
}
\end{remark}

\section{\'Equations de Monge-Amp\`ere et conjecture de Fujita}

Soit $X$ une variété projective de dimension $n$. La conjecture suivante a été formulée par T.~Fujita: \emph{si $L$ est un fibré ample sur $X$, alors le fibré adjoint
$K_X+ mL$ est engendre par ses sections globales et $K_X+ (m+1)L$ est très ample, dès que $m\geq n+1$.} Citons pour commencer le début de l'article \cite{ELN} de Ein-Lazarsfeld-Nakamaye, qui montre l'impact hors norme du travail \cite{Dem93} de Jean-Pierre dans cette direction : \emph{in the seminal paper [De1], Demailly drew on deep analytic tools to make the first serious attack on Fujita’s conjecture... While the numbers are rather far from Fujita’s predictions, this was the
first effective criterion for very ampleness, and it represented a real breakthrough.} Autrement dit, c'est Jean-Pierre qui a ouvert le bal ! Nous allons présenter maintenant quelques idées contenues dans l'article \cite{Dem93}.

\subsection{Le faisceau multiplicateur d'une fonction psh} Considérons un fibré en droites ample $L$ au-dessus d'une variété projective $X$. Le théorème \ref{DemGen} montre qu'afin de construire des sections holomorphes des fibrés type $K_X+ L$ qui ne s'annulent pas en un point $x_0\in X$, il suffit de montrer que $L$ admet une métrique singulière $h_L$ qui a les propriétés suivantes :
\begin{itemize}
\smallskip

\item le courant de courbure $T:= i\Theta(L, h_L)\geq 0$ est positif ;
\smallskip

\item le nombre de Lelong de $T$ en $x_0$ est supérieur à $n$ ;
\smallskip

\item il existe un voisinage $V\subset X$ de $x_0$ telle que le poids local $e^{-\varphi_L}$ de $h_L$ 
est localement intégrable en tout $y\in V\setminus \{x_0\}$. 
\end{itemize} 
\medskip

Afin de donner une condition suffisante pour que la condition d'intégrabilité du troisième point ci-dessus soit satisfaite, il convient de rappeler ici 
un concept fondamental qui a vu le jour dans \cite{Dem93}. C'est une manière d'associer à une fonction quasi-psh $\varphi$ un faisceau d'idéaux, c'est à dire un objet purement algébrique.

\begin{defn}\label{mult} Soit $\phi$ une fonction quasi-psh sur une variété hermitienne compacte $(X, \omega)$. L'ideal $\cI(\phi)\subset \mathcal O_X$ defini par
\[\cI(\phi)_x:= \{f\in \mathcal O_{X, x} : \int_{(X, x)}|f|^2e^{-\phi}d\lambda< \infty\}\]
s'appelle ideal multiplicateur de $\phi$. Si $L$ est un fibré en droites et $h_L$ est une métrique singulière sur $L$, dont le courant de courbure 
est minoré par $-C\omega$ ou $C\in \R$, l'ideal multiplicateur $\cI(h_L)$ de $h_L$ est défini via les poids locaux de cette métrique.
\end{defn}
\medskip

Certaines des propriétés importantes de l'idéal multiplicateur sont résumées dans l'énoncé suivant. 

\begin{thm}\label{DemMult}\cite{Dem93}
Soit $L$ un fibré holomorphe en droites sur une variété kählérienne compacte $(X, \omega)$, et soit $h_L$ une métrique singulière sur $L$
telle que $i\Theta(L, h_L)\geq -C\omega$. On a les affirmations suivantes.
\begin{enumerate}
\smallskip

\item[\rm (i)] L'idéal multiplicateur $\cI(h_L)$ est cohérent.
\smallskip

\item[\rm (ii)] Supposons que la constante $"C"$ ci-dessus est strictement négative. Alors le groupe de cohomologie 
$\displaystyle H^q\big(X, \mathcal O(K_X+ L)\otimes \cI(h_L)\big)=0$ s'annule, pour tout $q\geq 1$. 
\end{enumerate}
\end{thm}
\medskip

\begin{remark}{\rm 
Dans \cite{Dem93} cet énoncé est attribué à A. Nadel, cf. \cite{Nad}. Mais en regardant cette référence, on s'aperçoit que 
la définition de \emph{l'idéal multiplicateur} n'est pas vraiment celle de \ref{mult}. En effet, Nadel est intéressé par la construction des métriques de KE sur les variétés de Fano, et l'idéal qu'il introduit apparaît dans le contexte de la méthode de continuité dans les EDP, comme suit. \'Etant donnée une suite de fonctions $(\phi_k)_{k\geq 1}\subset \cC^{\infty}(X, \R)$ telle que 
\begin{equation}\label{eq11}\omega+ i\ddbar \phi_k\geq 0, \qquad \sup_X\phi_k= 0,\end{equation}
la fibre en $x\in X$ de l'idéal défini dans \cite{Nad} consiste en fonctions holomorphes $f\in \mathcal O_{X, x}$ telles que 
\[\sup_k \int_{(X, x)}|f|^2e^{-\phi_k}d\lambda< \infty.\]
C'est un peu déroutant : bien-entendu, on peux extraire une limite faible $\phi$ des $(\phi_k)_{k\geq 1}$ et considérer $\cI(\phi)$, compte tenu les hypothèses \eqref{eq11}, mais il n'est pas clair que les deux "idéaux multiplicateurs" coïncident. On se réfère a \cite{SiuDyn} pour quelques commentaires et résultats --potentiels-- au sujet de l'idéal multiplicateur de Nadel (ou il apparaît déguisé sous le nom de \emph{idéal multiplicateur dynamique}). 
} 
\end{remark}
%\smallskip

%\begin{remark}{\rm
%La notion de \emph{seuil log-canonique} d'un diviseur est fondamentale dans le programme des modèles minimaux : c'est comme cela qu'on localise les "pires" singularités d'un diviseur, ce qui permet de faire des raisonnements par recurrence très impressionnantes (voir par exemple la preuve du théorème de finitude de l'algèbre canonique d'une variété projective, \cite{}). La possibilité de formuler une notion analogue 
%pour les fonctions psh repose sur la solution de la conjecture d'ouverture par Berndtsson, cf. \cite{} et Guan-Zhou, cf. \cite{}. La conjecture a été 
%formulée par Demailly-Kollár dans \cite{}.
%}\end{remark}
\medskip 

Quoi qu'il en soit, sous les hypothèses du théorème précédent le morphisme 
\[H^0\big(X, \mathcal O(K_X+ L)\big)\to H^0\big(X, \mathcal O(K_X+ L)\otimes \mathcal O_X/\cI(h_L)\big)\]
est \emph{surjectif}, ce qui représente une amélioration du résultat \ref{DemGen}. Autrement dit, afin de produire une section $s$ 
de $K_X+L$ telle que $s_{x_0}\neq 0$, il suffit de montrer l'existence d'une métrique $h_L$ comme dans \ref{DemMult}, telle que 
$x_0$ soit un point isole parmi les zéros de $\cI(h_L)$. Dans cette direction, le résultat suivant dû à H. Skoda est très utile.

\begin{lemma}\label{SK}\cite{HS}
Soit $\varphi\in \Psh(\Omega)$ une fonction psh définie sur un ouvert de $\C^n$, et soit $x\in \Omega$ un point tel que 
$\nu(\varphi, x)< 1$. Alors la fonction $e^{-\varphi}$ est localement intégrable en $x_0$.  
\end{lemma}
\medskip

En conclusion, l'angle d'attaque de Jean-Pierre pour la conjecture de Fujita a été de construire une métrique $h_L$ dont les nombres de Lelong 
enregistrent un "saut" au voisinage d'un point donné $x_0\in X$. C'est à dire, si $T:= i/\pi\Theta(L, h_L)$, on souhaite que $\nu\big(T, x_0\big)$ soit grand (au moins $n$, afin que l'idéal multiplicateur soit non-trivial en ce point), et en même temps que $\nu\big(T, y\big)< 1$ pour tout $y\neq x_0$ dans un voisinage de $x_0$. Sans beaucoup de surprises, la première de ces deux conditions est beaucoup plus facile à satisfaire que la seconde, cf. la suite.

\subsection{Le théorème de Yau}  Le résultat fondamental suivant a été obtenu par S.-T. Yau dans \cite{Yau}.

\begin{thm}\label{Y}\cite{Yau} Soit $(X, \omega)$ une variété kählérienne compacte $n$-dimensionnelle, et soit $f\in \mathcal C^\infty(X, \R)$ une fonction reélle telle que 
\[\int_X\omega^n= \int_Xe^f\omega^n.\]
Alors il existe une fonction $\varphi \in \mathcal C^\infty(X, \R)$, unique à normalisation près, telle que
\[\omega+ i\ddbar \varphi> 0, \qquad (\omega+ i\ddbar\varphi)^n= e^f\omega^n.\]
\end{thm}
\medskip

Avant d'expliquer la manière dont cet énoncé est employé dans \cite{Dem93}\footnote{où on lui concède seulement le statut quasi-humiliant de "Lemma"...}, on peux pas résister à la tentation de faire un petit détour, comme suit.

\begin{remark} {\rm Un des aspects importants du théorème \ref{Y} est qu'il produit des métriques kählériennes dans une classe de 
cohomologie donnée qui reflètent les propriétés numériques du fibré canonique $K_X$ de $X$. Supposons par exemple que $-K_X$ soit nef. 
Dans \cite{DPS} on utilise \ref{Y} afin de montrer que cette hypothèse purement numérique a une contrepartie métrique :
il existe une famille de métriques kählériennes $(\omega_\ep)_{\ep> 0}$ dans la classe $\{\omega\}\in H^{1,1}(X, \R)$ 
définie par $\omega$ telle que $\Ricci_{\omega_\ep}\geq -\ep\omega_\ep$. Ainsi, on peux utiliser les méthodes de la géométrie des variétés riemanniennes à courbure de Ricci minorée afin d'investiguer les propriétés des variétés kählériennes dont le fibré anticanonique est nef, cf. \cite{MP1} et les références dedans. 
Une différence marquante par rapport au contexte riemannien est l'absence d'informations concernant le diamètre de $(X, \omega_\ep)$ lorsque $\ep\to 0$. En quelque sorte, cela est remplacé par la "normalisation" $\omega_\ep\in \{\omega\}$. Pour des développements importants du \cite{DPS} on se réfère au très beau travail de J. Cao et A. Höring cf. \cite{CH}.
}
\end{remark} 
\medskip

Considérons à présent un point $x_0\in X$ ainsi qu'un ouvert de coordonnées $x_0\in U\subset X$ qui contient $x_0$. Soient $z=(z_1,\dots, z_n)$ des coordonnées centrées en $x_0$. Un calcul direct basé sur le théorème de Stokes montre que 
\[\Big(\frac{i}{2\pi}\ddbar\log|z|^2\Big)^n= \delta_{x_0},\]
la distribution de Dirac en $x_0$. On voudrait utiliser $\delta_{x_0}$ à la place de l'élément volume $e^f\omega^n$ dans \ref{Y} (et espérer que la solution $\varphi$ sera au moins aussi singulière que $\log|z|^2$), mais on ne peux pas faire ceci directement.
La voie suivie dans \cite{Dem93} a été de construire d'abord une régularisation de $\delta_{x_0}$, comme suit.

Soit $\chi\in \mathcal C^{\infty}(\R, \R)$ une fonction convexe croissante, telle que $\chi(t)= t$ si $t\geq 0$ et $\chi(t)= -1/2$ si $t\leq -1$. 
On introduit la forme 
\[\alpha_\ep:= \frac{i}{2\pi}\ddbar\chi\big(\log|z|^2/\ep^2\big),\] 
semi-positive et définie sur $U$. \'Etant donné que $\alpha_\ep= \frac{i}{2\pi}\ddbar\log|z|^2)$ dès que $|z|> \ep$, le support de la puissance extérieure maximale $\alpha_\ep^n$ est contenu dans $(|z|\leq \ep)\Subset U$ ; on peux donc la considérer comme $(n,n)$--forme globale, $\mathcal C^\infty$ et semi-positive sur $X$, qui a les propriétés suivantes
\[\int_X\alpha_\ep^n= 1, \qquad \alpha_\ep^n \rightharpoonup \delta_{x_0}\]
(vérifiées en utilisant intégration par parties). Du coup, la famille $(\alpha_\ep^n)_{\ep>0}$ sera la régularisation de la distribution de Dirac mentionnée auparavant.

Dans le cadre de la conjecture de Fujita, on a un fibré ample $L\to X$, et soit $\omega_L\in c_1(X)$ une métrique kählérienne.
Soit $\tau> 0$ un nombre réel arbitraire, tel que 
\[\tau^n< \int_X\omega_L^n = L^n. \]
D'après \ref{Y}, il existe une unique fonction $\varphi_{\ep}\in \mathcal C^\infty(X, \R)$ telle que
\begin{equation}\label{eq16}\sup_X\varphi_\ep= 0, \qquad \omega_L+ \frac{i}{2\pi}\ddbar \varphi_\ep>0\end{equation}
et telle que l'égalité
\begin{equation}\label{eq15}
(\omega_L+ \frac{i}{2\pi}\ddbar \varphi_\ep)^n= \tau^n\alpha_\ep^n+ \big(1-\frac{\tau^n}{L^n}\big)\omega_L^n
\end{equation}
entre les deux éléments de volume soit satisfaite, pour tout $\ep> 0$. 

Maintenant, la normalisation et la condition de positivité dans \eqref{eq16} font qu'on peux extraire une limite faible
\begin{equation}\label{eq17}
T= \lim_{\ep\to 0}\omega_\ep, \qquad T\in c_1(L).
\end{equation}
où on note $\displaystyle \omega_\ep:= \omega_L+ \frac{i}{2\pi}\ddbar \varphi_\ep$ la solution de l'équation dans le théorème \ref{Y}.
On esquisse maintenant l'argument montrant que 
\begin{equation}\label{eq18}
\varphi_\ep\leq \tau \log(|z|^2+ \ep^2)
\end{equation}
en suivant \cite{Dem93}. D'abord, on peux supposer que $U$ est biholomorphe à la boule unité dans $\C^n$ et que $\omega_L|_U= \frac{i}{2\pi}\ddbar \phi_U$ quitte à restreindre $U$, où $\phi_U$ est une fonction de classe $\cC^\infty$ définie sur un voisinage de $\ol U$. L'équation \eqref{eq15} montre qu'on a
\begin{equation}\label{eq19}
\big(i\ddbar (\phi_U+ \varphi_\ep)\big)^n> \big(\tau i \ddbar\chi\big(\log|z|^2/\ep^2\big)\big)^n
\end{equation}
en chaque point de $U$. Par ailleurs, il existe une constante $C> 0$ telle que pour tout point $z$ situé au bord $\partial U$ de $U$ l'inégalité
\begin{equation}\label{eq20}
\phi_U(z)+ \varphi_\ep(z)\leq \tau \big(\chi\big(\log|z|^2/\ep^2\big)+ \log \ep^2\big)+ C
\end{equation}
soit vérifiée, dès que $0< \ep\ll 1$. Ceci est clair, car la restriction de $\chi$ à l'ensemble des réels positifs est la fonction identité. 
En utilisant \eqref{eq19}, on déduit que le maximum de la fonction
\[z\to \phi_U(z)+ \varphi_\ep(z)- \tau \big(\chi\big(\log|z|^2/\ep^2\big)+ \log \ep^2\big)- C, \qquad z\in \ol U\] 
ne peux pas être atteint à l'intérieur de $U$. 

En conjonction avec \eqref{eq20} cela montre que 
$$\phi_U(z)+ \varphi_\ep(z)\leq \tau \big(\chi\big(\log|z|^2/\ep^2\big)+ \log \ep^2\big)+ C$$ pour tout point $z\in U$, et donc on obtient
\begin{equation}\label{eq21}
\phi_U(z)+ \varphi_\ep(z)\leq \tau\log|z|^2+ C
\end{equation}
pour tout $|z|> \ep$. L'inégalité \eqref{eq18} s'en suit, et on a montré l'énoncé suivant.

\begin{cor}\label{K1}
Soit $L\to X$ un fibré en droites ample, et soit $\tau\in \R_+$ tel que $\tau^n< L^n$. Pour tout point $x_0\in X$ 
il existe une métrique singulière $h_L$ de $L$ telle que :
\begin{enumerate}
\smallskip

\item[\rm (a)] le courant de courbure $i\Theta(L, h_L)\geq 0$ est positif,
\smallskip

\item[\rm (b)] si on note $\phi$ le poids de $h_L$ au voisinage de $x_0$, on a $\phi(z)\leq \tau\log|z|^2+ \mathcal O(1)$,
où $z$ sont des coordonnées centrées locales en $x_0$.  
\end{enumerate}
\end{cor}
\smallskip

\begin{remark}\label{R1}{\rm 
Le corollaire précédent montre qu'il est possible de construire des métriques singulières sur $L$ avec un pôle logarithmique
en un point fixé à l'avance. Mais ces métriques sont obtenues de manière "abstraite", en considérant des limites faibles
et ainsi, il semble très difficile de contrôler leurs singularités aux points $y\neq x_0$. Nous allons expliquer en quelques mots 
la manière de procéder dans \cite{Dem93} dans la section suivante.}
\end{remark}

\subsection{Une inégalité d'auto-intersection} Supposons à présent que l'auto-intersection du fibré ample $L$ vérifie la condition suivante
\begin{equation}\label{eq123}L^n> n^n.\end{equation}
Si $x_0\in X$ est un point fixé, il existe une métrique singulière $h_L$ sur $L$, dont le nombre de Lelong du courant de courbure $T= i\Theta(L, h_L)$ en $x_0$ est supérieur à $n$, d'après le corollaire \ref{K1}. Cela implique en particulier que la fibre de l'idéal multiplicateur $\cI(h_l)_{x_0}$ est non-triviale, i.e. strictement contenue dans $\O_{X, x_0}$. Comme on l'a déjà évoqué dans la remarque précédente \ref{R1}, a priori il semble impossible de contrôler $\nu(T, y)$ pour $y\neq x_0$.

Considérons la décomposition de Siu du courant $T$, i.e.
\[T= \sum_{i\geq 1} \nu_i[Y_i]+ R.\]
Il peux très bien arriver que $x_0$ se trouve sur une (ou plusieurs) des hypersurfaces $Y_i$ telles que $\nu_i> n$ et dans ce cas, on aura $\nu(T, y)> n$ pour tout $y\in Y_i$. L'énoncé suivant est un premier pas vers la solution de ce problème.

\begin{lemma}\label{L1}\cite{Dem93} Sous les hypothèses du corollaire \ref{K1},
  on a \[\int_X R\wedge \omega_L^{n-1}\geq \big(1-\frac{\tau^n}{L^n}\big)^{\frac{1}{n}}\int_X\omega_L^n\]
  lorsque $T= \lim_\ep\omega_\ep$ est une limite faible extraite de $\displaystyle (\omega_\ep)_{\ep> 0}$, cf. équation \eqref{eq15}.
\end{lemma}
Avant d'esquisser l'argument montrant cette inégalité,
remarquons que cela implique
\[\sum \nu_iL^{n-1}\cdot Y_i\leq \Big(1-\big(1-\frac{n^n}{L^n}\big)^{\frac{1}{n}}\Big)L^n,\]
et si on suppose de plus que $L$ est suffisamment positif tel qu'on ait
\begin{equation}\label{eq22}
L^{n-1}\cdot Y\geq \Big(1-\Big(1-\frac{n^n}{L^n}\Big)^{\frac{1}{n}}\Big)L^n 
\end{equation}
pour toute hypersurface $Y\subset X$, alors $\nu_i< 1$. Donc, on a montré le corollaire suivant.
\begin{cor}\label{K2}
  Supposons que $L$ est un fibré en droites ample tel que les inégalités \eqref{eq123} et \eqref{eq22} soient satisfaites. Pour tout point $x_0\in X$ il existe une métrique singulière $h_L$ de $L$ telle que 
  \begin{itemize}
  \item le courant de courbure $T$ est semi-positif et $x_0\in E_n(T)$ ; 
  
  \item la codimension de $E_1(T)$ dans $X$ est au moins deux. \end{itemize}
\end{cor}

Par exemple, dans le cas des surfaces cela signifie que $x_0$ est un point isolé dans $E_1(T)$ et donc le fibré adjoint $K_X+ L$ est globalement engendré dès que les conditions numériques dans \ref{K2} sont satisfaites.
\medskip

Pour ce qui est du lemme \ref{L1}, remarquons qu'on a l'inégalité
\begin{equation}\label{eq23}\int_X\theta T\wedge \omega_L^{n-1}\geq \lim\inf_{\ep\to 0}\int_X\theta\omega_\ep\wedge \omega_L^{n-1}\end{equation}
pour toute fonction positive $\theta$, compte tenu du fait que $T$ est une limite faible des $(\omega_\ep)_{\ep>0}$.
Par ailleurs, on a $\displaystyle \omega_\ep\wedge \omega_L^{n-1}= \frac{1}{n}\tr_{\omega_L}(\omega_\ep)\omega_L^n\geq \Big(\frac{\omega_\ep^n}{\omega_L^n}\Big)^{1/n}\omega_L^n
\geq \big(1-\frac{\tau^n}{L^n}\big)^{\frac{1}{n}}\omega_L^n$, qu'on obtient facilement par l'inégalité de la moyenne combinée avec l'équation \eqref{eq15}. Donc, l'inégalité \eqref{eq23} implique
\begin{equation}\label{eq24}\int_X\theta T\wedge \omega_L^{n-1}\geq
  \big(1-\frac{\tau^n}{L^n}\big)^{\frac{1}{n}} \int_X\theta \omega_L^{n}\end{equation}
et le lemme \ref{L1} s'en suit par un choix adéquate de la fonction $\theta$.
\medskip

Malheureusement, l'analyse des singularités de $T$ le long des ensembles de codimension supérieure à un est beaucoup plus laborieuse. La raison principale réside dans le fait qu'en général, il est impossible de définir le produit d'intersection  $T\wedge T$ en tant que courant positif de bidegré $(2,2)$ (cela reviendrait à multiplier des mesures...). L'observation importante est que si 
\[S= \beta+ i\ddbar \psi\geq 0\]
est un courant positif fermé dans la classe $\{\beta\}$ tel que $\psi$ est localement borné sur $X\setminus Z$, avec $\codim_X Z\geq 2$, alors les techniques de Bedford-Taylor perfectionnées par Jean-Pierre permettent de définir le produit extérieur 
$T\wedge S$ : c'est un courant positif fermé de bidegré $(2,2)$, et sa classe de cohomologie coïncide avec $\{T\}\wedge \{S\}$. Concernant les nombres de Lelong de $T\wedge S$, on montre dans \cite{Dem93} que 
\begin{equation}\label{eq25}
T\wedge S\geq \sum_k \nu_k\mu_k[Z_k]
\end{equation}
où les ensembles analytiques $Z_k\subset Z$ sont les composantes de $Z$ de codimension précisément deux, et $\nu(T, x)= :\nu_k, \nu(S, x)= :\mu_k$ sont les nombres de Lelong de $T$ et respectivement $S$ au point générique $x\in Z_k$. L'inégalité \eqref{eq25} peut-être interprétée comme une version très généralisée du théorème classique de Bézout.

Du coup, même si on peux pas construire directement le produit d'intersection $T\wedge T$, il est envisageable de "se débarrasser" en préalable des singularités de $T$ en codimension un, et utiliser ce nouveau courant dans \eqref{eq25}. C'est dans cet esprit  que Jean-Pierre a démontré le résultat fondamental suivant.

\begin{thm}\label{DemReg}\cite{DemReg}
Soit $(X, \omega)$ une variété complexe compacte munie d'une métrique hermitienne,  
et soit $T= \alpha+ i\ddbar \varphi$ un courant positif fermé. Pour tout réel $c> 0$ et tout $k\geq 1$ 
il existe une fonction quasi-psi $\varphi_{c, k}\in \mathcal C^\infty\big(X\setminus E_c(T)\big)$ avec les propriétés suivantes. 
\begin{enumerate}
\smallskip

\item[\rm (a)] Le courant $T_{c, k}:= \alpha+ i\ddbar \varphi_{c, k}$ est plus grand que $-\ep_k\omega- c \tau_X$, ou $\ep_k\to 0_+$
et $\tau_X\geq 0$ est une $(1,1)$-forme semi-positive qui dépend de la géométrie de $X$ seulement.
\smallskip

\item[\rm (b)] En chaque point $x\in X$ on a $\displaystyle \nu(T_{c, k}, x)= \max\big(\nu(T, x)-c, 0\big)$.
\end{enumerate}
\end{thm} 

L'énoncé dans \cite{DemReg} est beaucoup plus complet, mais même cette version édulcorée permet de voir que si $\displaystyle 0<c\leq \max_{x\in X}\nu(T, x)$ est choisi tel que les composantes irréductibles de $E_c(T)$ ont codimension au moins deux, alors 
le produit d'intersection $T\wedge (T_{c, k}+ \ep_k\omega+ c \tau_X)$ est bien défini, et de plus on a une borne inférieure pour les coefficients qui apparaissent dans la partie "cycle en codimension deux" de la decomposition de Siu de ce courant. 

Un coup d'oeil attentif à la manière dont le cas des composantes de codimension un a été traité montrera que l'inégalité obtenue dans le lemme \ref{L1} est très importante : sans cela, on a aucune chance de montrer que $\nu_k< 1$, même si on remplace le fibré $L$ par une de ses puissances. Afin de formuler la généralisation de cette inégalité, remarquons d'abord que tout courant positif fermé admet une décomposition de Lebesgue 
\[T= T_{\rm s}+ T_{\rm abc}\]
via ses coefficients (dans les notations \eqref{eq4}, on obtient cela en considérant la décomposition de Lebesgue de 
chaque mesure $\mu_{j\ol k}$). Notons qu'en général, les deux composantes de $T$ ci-dessus ne sont plus fermées ; néanmoins, sous les hypothèses du lemme \ref{L1} on a l'inégalité 
\[\int_X T_{\rm abc} \wedge \omega_L^{n-1}\geq \big(1-\frac{\tau^n}{L^n}\big)^{\frac{1}{n}}\int_X\omega_L^n.\]
 \smallskip
 
 On introduit à présent les invariants $b_p:= \inf\{c> 0 : \codim \big(E_c(T), x\big)\geq p, \forall x\in X\}$, ou $p=1,\dots, n+1$ ; l'énoncé qui permet d'obtenir des bornes effectives $m_n$ pour l'existence des sections $K_X+ m_nL$ qui ne s'annulent pas en un point donné est une version généralisée du théorème suivant (on reprend les notations du début de cette section). 
 
 \begin{thm}\cite{Dem93} Soient $(Z_k)_{k\geq 1}$ les composantes irréductibles de codimension deux contenues dans $\bigcup_{c> b_2}E_c(T)$ ; on note $\nu_k$ le nombre de Lelong de $T$ au point générique de $Z_k$. Alors il existe un courant positif $\Theta_2\in \{T\}\wedge \{T+ b_2\tau_X\}$ de bidegré $(2,2)$ tel que l'inégalité
\begin{equation}\label{eq26}
\Theta_2\geq \sum_k\nu_k(\nu_k- b_2)[Z_k]+ T_{\rm abc}\wedge (T_{\rm abc}+ b_2\tau_X)
\end{equation} 
est satisfaite au sens des courants.
 \end{thm}  
 
Pour finir cette section, mentionnons que dans \cite{Dem93}, \cite{DemReg} 
on construit pour chaque $p\geq 2$ un courant positif
\[\Theta_p\in \{T\}\wedge \{T+ b_2\tau_X\}\wedge\dots\wedge\{T+ b_p\tau_X\}\] vérifiant une inégalité similaire au \eqref{eq26}. 
Le lecteur trouvera dans \cite{Dem93} une brillante manière de combiner tout cela --et de se débarrasser de la forme $\tau_X$, qui est "hors jeu", pour ainsi dire-- afin d'aboutir au résultat principal. Cet article contient une abondance d'idées et techniques très intéressantes, tout cela est franchement beau!\footnote{et quelque part injuste, car comme on a déjà mentionné, les bornes effectives obtenues sont loin des predictions de Fujita...} 

\subsection{Développements ultérieurs} L'article de Jean-Pierre a attiré l'attention de la \emph{fine lame de la géométrie algébrique}, pour s'exprimer ainsi. Autre que l'article \cite{ELN} déjà mentionné\footnote{et qui se trouve à l'origine des recherches concernant les invariants asymptotiques des séries linéaires}, citons les contributions fondamentales de Ein-Lazarsfeld \cite{EL93}, J.~Kollár \cite{Kol93}, Y.~Kawamata \cite{Kaw97} et Angehrn-Siu \cite{AS95}. 
\`A ce jour, la question de déterminer une constante effective $m_n$ dépendant de la dimension $n= \dim(X)$ de la variété $X$ seulement 
telle que $K_X+ m_nL$ est très ample dès que $L$ est ample reste à élucider. 

\section{Le c\^one de K\"ahler}
Soit $X$ une variété projective et soit $L\to X$ un fibré en droites au-dessus $X$. Le théorème suivant, le très important critère d'amplitude de Nakai-Moishezon s'énonce comme suit : \emph{$L$ est ample si et seulement si les nombres d'intersection 
$\displaystyle L^d\cdot Y= \int_{Y_{\rm reg}}c_1(L)^d> 0$ de $L$ avec toute sous-variété réduite, irréductible $Y\subset X$ de dimension $d$ est strictement positif.} Probablement motivé par cela, A. Beauville a formulé la conjecture suivante :  
  
\begin{conjecture}\label{Beau}
Soit $X$ une variété projective, et soit $\{\alpha\}\in H^{1,1}(X, \R)$ une classe réelle de type $(1,1)$. Alors $\{\alpha\}$ est une classe 
kählérienne si et seulement si on a $\displaystyle \int_{Y_{\rm reg}}\alpha^d> 0$ pour toute sous-variété réduite, irréductible $Y\subset X$ de dimension $d$.
\end{conjecture}  
\medskip

La version kählérienne de cette conjecture a été obtenue dans \cite{DP04}, comme suit.
\begin{thm}\label{NakMoi} %\cite{DP04}
Soit $(X, \omega)$ une variété kählérienne compacte. Alors son cône de Kähler $\mathcal K\subset H^{1,1}(X, \R)$ est une des composantes connexes de l'ensemble
\begin{equation}\label{eq30}\Big\{\{\alpha\}\in H^{1,1}(X, \R) : \int_{Y_{\rm reg}}\alpha^d> 0, \forall Y\subset X, \dim Y= d\Big\}.\end{equation}
De plus, on a l'égalité
\begin{equation}\label{eq127}
\mathcal K= \Big\{\{\alpha\}\in H^{1,1}(X, \R) : \int_{Y_{\rm reg}}\alpha^k\wedge \omega^{d-k}> 0, \forall Y\subset X, \dim Y= d, k=0,\dots, d\Big\}
\end{equation}
\end{thm}
\medskip

Notons que la présence de la métrique de référence $\omega$ dans \eqref{eq127} ne signifie pas que $\mathcal K$ dépend de cette métrique : son rôle est de "détecter" le cône de Kähler parmi les composantes de \eqref{eq30}, au cas ou $X$ ne possède pas assez de sous-ensembles analytiques (voir e.g. le cas d'un tore de dimension algébrique zéro).
Nous allons présenter maintenant les étapes principales de la preuve de ce théorème. 

\subsection{Concentration de la masse et méthode de continuité}
Nous rappelons ici les notions suivantes, qui représentent les contre-parties "transcendantes" des définitions usuelles en géométrie algébrique. 

\begin{defn}\label{nef} Soit $\{\alpha\}\in H^{1,1}(X, \R)$ une classe réelle de type $(1,1)$ sur une variété complexe compacte $(X, \omega)$.
\begin{enumerate}
\smallskip

\item[\rm (1)] On dit que $\{\alpha\}$ est nef s'il existe une famille $(f_\ep)_{\ep> 0}\subset \cC^\infty(X, \R)$ de fonctions de classe $\cC^\infty$ telle que $\displaystyle \alpha+ i\ddbar f_\ep\geq -\ep\omega$ pour tout $\ep> 0$.
\smallskip

\item[\rm (2)] La classe $\{\alpha\}$ est nef et big si $\{\alpha\}$ est nef, et de plus $\displaystyle \int_X\alpha^n> 0$.
\end{enumerate}
\end{defn}

Comme on a déjà vu dans le théorème \ref{DemDict}, si $X$ est projective et si $\{\alpha\}= c_1(L)$, alors cette classe est nef au sens de la définition ci-dessus si et seulement si $L$ est nef, cf. définition \ref{mo}. Si $X$ est une variété kählérienne compacte, alors la classe $\{\alpha\}$ est nef si
et seulement si $\{\alpha\}\in \ol{\mathcal K}$, i.e. le cône nef est l'adhérence du cône de Kähler de $X$ dans $H^{1,1}(X, \R)$.
\medskip

Dans ce contexte, on a le résultat suivant.
\begin{thm}\label{NefBig}\cite{DP04}
Soit $(X, \omega)$ une variété kählérienne compacte, et soit $\{\alpha\}\in H^{1,1}(X, \R)$ une classe nef et big. Considérons également un 
sous-ensemble analytique irréductible $Y\subset Y$ tel que $\codim_XY= p$. Alors la classe $\{\alpha^p\}$ contient un courant positif fermé
$\Theta_Y$ tel que l'inégalité suivante 
\[\Theta_Y\geq \ep_0[Y]\]
soit satisfaite pour un certain réel positif $\ep_0> 0$. 
\end{thm}

\`A la première vue, l'énoncé ci-dessus \ref{NefBig} n'est pas entièrement satisfaisant, compte tenu le mélange des hypothèses numériques et métriques. Cependant, le théorème \ref{NakMoi} est une conséquence presque immédiate de \ref{NefBig}, grâce à l'observation suivante.
On peut supposer que $\alpha+ \omega$ est positive (quitte à remplacer $\omega$ par un multiple assez grand). Considérons l'ensemble
\[I= \big\{t\in [0,1] : \{\alpha+ t\omega\} \in \mathcal K\big\} ;\]
il est clair que $I\neq\emptyset$ car $1\in I$, et est également que $I$ est un ensemble ouvert. On doit montrer que $I$ est fermé, et soit $t_\infty\in \ol I$. 

La classe limite $\{\alpha + t_\infty\omega\}$ se trouve --au pire-- dans l'adhérence du cône de Kähler de $X$, 
autrement dit c'est une classe \emph{nef}. L'auto-intersection de celle classe est strictement positive, car c'est le cas du premier terme de la somme de droite
\[\int_X(\alpha+ t_\infty\omega)^n= \int_X\alpha^n+ \sum_{k=1}^nt_\infty^{n-k}{{n}\choose {k}}\int_X\alpha^k\wedge \omega^{n-k},\]
et les autres sont positifs ou nuls, (conséquence du $t_\infty\geq 0$ et \eqref{eq127}). 

Considérons à présent la variété produit $X\times X$, avec les projections $\pi_i:X\times X\to X$, ou $i=1, 2$ sur les deux facteurs.
La classe 
\[\beta:= \pi_1^\star(\alpha+ t_\infty\omega)+ \pi_2^\star(\alpha+ t_\infty\omega)\]
est également nef et big, et le théorème \ref{NefBig} montre l'existence d'un courant $\Theta\in \{\beta\}^n$ tel que
\[\Theta\geq \ep_0[\Delta],\]
où on désigne par $[\Delta]$ le courant d'intégration sur la diagonale $\Delta\subset X\times X$. 

Soit $\displaystyle \delta:= n\int_X(\alpha+ t_\infty\omega)^{n-1}\wedge\omega$ ; c'est un nombre réel strictement positif, et il est facile de voir que le 
courant image directe $\displaystyle T:= \frac{1}{\delta}\pi_{1\star}(\Theta\wedge\pi_2^\star\omega)$ a les propriétés suivantes
\[T\in \{\alpha+ t_\infty\omega\}, \qquad T\geq \ep_1\omega,\] 
ou $\ep_1>0$. 

En résumé, on a montré l'existence d'un courant kählérien $T$ dans la classe qui nous intéresse, ce qui représente le pas décisif vers la fin de la preuve. En effet, si $X$
est de dimension algébrique zéro (i.e. s'il n'existe pas des sous-ensembles analytiques de dimension positive de $X$), alors le théorème de régularisation \ref{DemReg} montre immédiatement que la classe $\{\alpha+ t_\infty\omega\}$ est kählérienne. En général, on peux 
remplacer $T$ par un courant kählérien dont les singularités sont concentrés le long d'un ensemble analytique de $X$, et conclure par récurrence
.\footnote{Bon, à dire vrai l'histoire est un peu plus longuette que cela, on se réfère à \cite{DP04} pour les détails.}
\medskip

Pour ce qui est de la preuve du théorème \ref{NefBig}, on utilise dans \cite{DP04} une fois de plus l'équation de Monge-Ampère. Considérons un recouvrement de $X$ par des ouverts de coordonnées $(U_j)_{j=1,\dots, N}$. Pour chaque $j$, soit $(g_{jk})_k\subset \mathcal O(U_j)$ un système de générateurs de l'idéal correspondant a $Y$. \'Etant donné $\ep>0$, on introduit la fonction 
\[\psi_\ep: X\to \R, \qquad \psi_\ep(z)= \log\Big(\ep^2+ \sum_{j,k}\theta_j(z)^2|g_{jk}(z)|^2\Big),\]
où $(\theta_j)_{j=1,\dots, N}$ est une partition de l'unité, telle que $\Supp(\theta_j)\subset U_j$ pour chaque indice $j$. 
Il existe une constante $C> 0$ telle que l'inégalité
\[\rho_\ep:= \omega+ \frac{1}{C}i\ddbar\psi_\ep \geq \frac{1}{2}\omega\]
soit vérifiée sur $X$, pour tout valeur du paramètre $\ep$. Cette affirmation est prouvée par un calcul direct, cf. \cite{DP04}, et cela représente une version "régularisée" d'un résultat général du à H.~Grauert, cf. \cite{HG}.

Soit $f_\ep\in\mathcal C^\infty(X, \R)$ la fonction définie par la relation suivante
\[\rho_\ep^n= e^{f_\ep}\omega^n.\]
Pour chaque $\ep> 0$ considérons la constante $c_\ep:= \frac{\int_X(\alpha+ \ep\omega)^n}{\int_X\omega^n}$.
L'hypothèse \emph{$\alpha$ nef et big} montre l'existence d'un réel $C> 1$ tel que
\[C^{-1}< c_\ep< C\]   
pour tout $\ep> 0$.
Le théorème \ref{Y} montre l'existence d'une fonction $\varphi_\ep\in\mathcal C^\infty(X, \R)$ telle que les relations suivantes
\begin{equation}\label{eq27}
\alpha+ \ep\omega+ i\ddbar \varphi_\ep> 0, \qquad (\alpha+ \ep\omega+ i\ddbar \varphi_\ep)^n= c_\ep e^{f_\ep}\omega^n, \qquad 
\end{equation}
soient satisfaites, pour tout $\ep> 0$. 

La dernière étape de la preuve consiste à considérer une limite faible 
\begin{equation}\label{eq28}
\Theta_Y:= \lim_{\ep\to 0}(\alpha+ \ep\omega+ i\ddbar \varphi_\ep)^p,
\end{equation}
et à montrer que $\Theta_Y$ est supérieur à un multiple positif du courant d'intégration sur $Y$. Nous avons déjà discuté un résultat similaire 
dans le cas d'un point dans le section précédente. Il va de soi que traiter le cas $\dim Y> 0$ est plus difficile mais somme toute, l'idée importante est l'utilisation de l'équation de MA \eqref{eq27} ; on renvois à \emph{loc. cit.} pour les détails (qui consistent a montrer que le produit de $\Theta_Y$
avec la fonction caractéristique de $Y$ est non-nul).
\medskip

Une conséquence intéressante du théorème \ref{NakMoi} c'est le comportement du cône de Kähler par déformation.

\begin{thm}\cite{DP04}
Soit $\mathcal X\to \mathcal S$ une famille de variétés kählériennes compactes sur une base irréductible $\mathcal S$. Les cônes de Kähler
$\mathcal K_t\subset H^{1,1}(\mathcal X_t, \C)$ sont invariants par transport parallèle par rapport à la composante $\nabla^{1,1}$ de type $(1,1)$ de la connexion de Gau\ss-Manin pour $t\in \mathcal S$ très générique (dans le complément d'une réunion dénombrable d'ensembles analytiques propres de $\mathcal S$).
\end{thm}

\subsection{Développements ultérieurs} Une très jolie preuve alternative du résultat principal dans \cite{DP04} a été trouvé par I. Chiose dans \cite{Jose}. Par ailleurs, les résultats établis dans \cite{DP04} ont motivé l'étude des propriétés géométriques des classes réelles de type $(1,1)$, par analogie avec les fibrés en droites. C'est un domaine qui reste encore très actif de nos jours, voici quelques \emph{classiques} dans cette direction \cite{Bou}, \cite{CT22}, \cite{DWN}. Aussi, il se trouve que ce type de résultats sont indispensables dans les développements récents de la version kählérienne du MMP, \cite{HP}, \cite{DHP}.

\end{document}